\documentclass[12pt]{article}

\usepackage{amsmath}

\usepackage{amssymb}

\usepackage{amsbsy}

\begin{document}

\baselineskip=20pt

\date{}

\title{{\sc The Div-Curl Lemma Revisited}}

\author{{\sc Dan Poli\v sevski $^\dag$}}

\maketitle

\vspace{0.5cm}

{\bf Abstract.} The Div-Curl Lemma, which is the basic result of the
compensated compactness theory in Sobolev spaces(see [2]-[6]), was
introduced by [1] with distinct proofs for the $L^2(\Omega)$ and
$L^p(\Omega), p \neq 2,$ cases. In this note we present a slightly
different proof, relying only on a Green-Gauss integral formula and
on the usual Rellich-Kondrachov compactness properties.



\vspace{0.5cm}

\section{An integral formula}

From now on, for any vector distribution $w$ in $\mathbb{R}^N$, $N
\geq 2$, we denote:

\vspace{0.3cm}

\noindent$\mbox{D}(w) = \sum\limits_{i = 1}^n \;\dfrac{\partial
w_i}{\partial x_i}$\;\;\;\; and \;\;\;\;$\mbox{C}_{ij}(w) =
\dfrac{\partial w_i}{\partial x_j} - \dfrac{\partial w_j}{\partial
x_i}$, $1 \leq i, j \leq N$.

\vspace{0.3cm}

{\bf Theorem 1.1.} {\it Let $\Omega$ be an open bounded set of
$\mathbb{R}^N$, $N \geq 2$. Let $p, q > 1$ such that $1/p + 1/q =
1$. If $u \in [W_{loc}^{2, p}(\Omega)]^N$ and $v \in [W_{loc}^{2,
q}(\Omega)]^N$, then for any $\varphi \in \mathcal{D}(\Omega)$ we
have:

$$
\begin{array}{l}
\int\limits_\Omega \varphi \bigtriangleup u \bigtriangleup v  = -
\langle \; \mbox{D}(\bigtriangleup u), \varphi \; \mbox{D}(v)\;
\rangle_p - \int\limits_{\Omega}
\;\mbox{D}(u)\; (\bigtriangledown \varphi) \bigtriangleup v \;+ \\[7mm]
+ \int\limits_{\Omega} \; \mbox{D}(u) \;(\bigtriangledown
\varphi)\bigtriangledown \mbox{D}(v) - \int\limits_{\Omega}
\; \mbox{D}(v) \;(\bigtriangledown \varphi)\bigtriangledown \mbox{D}(u) \; - \\[5mm]
- \sum\limits_{1 \leq i, j \leq N} \; \int\limits_{\Omega}
\;\bigtriangleup v_i \;\mbox{C}_{ij}(u)\;\dfrac{\partial
\varphi}{\partial x_j} - \dfrac{1}{2} \sum\limits_{1 \leq i, j \leq
N} \; \langle \;\mbox{C}_{ij}(\bigtriangleup v), \varphi\;
\mbox{C}_{ij}(u)\; \rangle_q ,
\end{array} \leqno(1.1)
$$
where $ \; \langle \; , \; \rangle_p \;$  denotes the pairing of
$W^{-1, p}(\Omega)$ with $W_0^{1, q} (\Omega)$ and $ \; \langle \; ,
\; \rangle_q \;$ denotes the pairing of $W^{-1, q}(\Omega)$ with
$W_0^{1, p}(\Omega)$.}

\vspace{0.5cm}

{\bf \emph{Proof.}} For any $u \in [W_{loc}^{2, p}(\Omega)]^N, v \in
[W_{loc}^{2, q}(\Omega)]^N$ and $\varphi \in \mathcal{D}(\Omega)$ we
have
$$
\begin{array}{l} \langle \;\mbox{D}(\bigtriangleup u), \varphi \;
\mbox{D}(v)\; \rangle_p + \dfrac{1}{2} \sum\limits_{1 \leq i, j \leq
N} \; \langle \;\mbox{C}_{ij}(\bigtriangleup v),
\varphi \; \mbox{C}_{ij}(u)\; \rangle_q =\\[8mm] = -
\int\limits_{\Omega} \; \bigtriangleup u \bigtriangledown (
\;\varphi \; \mbox{D}(v)\;) + \int\limits_\Omega \;\bigtriangleup v
\;(\;\bigtriangledown \mbox{D}(u) - \bigtriangleup u\;)\; \varphi \; - \\[5mm]
- \sum\limits_{1 \leq i, j \leq N}\; \int\limits_\Omega
\;\bigtriangleup v_i \;\mbox{C}_{ij}(u)\;\dfrac{\partial
\varphi}{\partial x_j}\;, \end{array} \leqno (1.2)
$$
which implies:
$$
\begin{array}{l}
\int\limits_\Omega \varphi \bigtriangleup u \bigtriangleup v  = -
\langle \; \mbox{D}(\bigtriangleup u), \varphi \; \mbox{D}(v)\;
\rangle_p - \int\limits_{\Omega}
\;\mbox{D}(u)\;(\bigtriangledown \varphi)\bigtriangleup v\;- \\[7mm]
- \int\limits_{\Omega} \; \bigtriangleup u \bigtriangledown (
\;\varphi \; \mbox{D}(v)\;)
 + \int\limits_{\Omega} \;\bigtriangleup v \bigtriangledown
(\;\varphi \;\mbox{D}(u)\;) \; - \\[5mm]
- \sum\limits_{1 \leq i, j \leq N} \; \int\limits_{\Omega}
\;\bigtriangleup v_i \;\mbox{C}_{ij}(u) \;\dfrac{\partial
\varphi}{\partial x_j} - \dfrac{1}{2} \sum\limits_{1 \leq i, j \leq
N} \; \langle \;\mbox{C}_{ij}(\bigtriangleup v), \varphi\;
\mbox{C}_{ij}(u)\;\rangle_q .
\end{array} \leqno(1.3)
$$

As $(\bigtriangledown \mbox{D}(u) - \bigtriangleup u)\in
[L^p_{loc}(\Omega)]^N$ is divergence-free and as $\varphi
\mbox{D}(v)\in W_0^{1, q} (\Omega)$, it follows:
$$
\begin{array}{l} \int\limits_{\Omega}  \bigtriangleup u
\bigtriangledown ( \varphi \; \mbox{D}(v))= - \int\limits_{\Omega}
\; (\bigtriangledown \; \mbox{D}(u) - \bigtriangleup u)
\bigtriangledown ( \varphi \; \mbox{D}(v))
 + \\[5mm]+\int\limits_{\Omega} \;
\bigtriangledown \; \mbox{D}(u) \bigtriangledown (\varphi
\;\mbox{D}(v))= \int\limits_{\Omega} \; \bigtriangledown \;
\mbox{D}(u) \bigtriangledown (\varphi \;\mbox{D}(v)).
\end{array}\leqno(1.4)
$$
By treating similarly the next term of (1.3), that is
$\;\int\limits_{\Omega}\bigtriangleup v \bigtriangledown (\varphi
\mbox{D}(u))\;$, we get straightly (1.1).

\vspace{0.5cm}

\section{The Div-Curl Lemma}

{\bf Theorem 2.1.} {\it Let $\Omega$ be an open bounded set of
$\mathbb{R}^N$, $N \geq 2$. Let $p, q > 1$ such that $1/p + 1/q =
1$. For any $\varepsilon > 0$, let $a^\varepsilon \in
[L^p(\Omega)]^N$ and $b^\varepsilon \in [L^q(\Omega)]^N$ with the
properties:

\vspace{0.1cm}

(i) $a^\varepsilon \rightharpoonup a$ weakly in $[L^p(\Omega)]^N$,
as $\varepsilon \to 0$,

\vspace{0.1cm}

(ii) $b^\varepsilon \rightharpoonup b$ weakly in
$[L^q(\Omega)]^N$, as $\varepsilon \to 0$,

\vspace{0.1cm}

(iii) $\{\mbox{D}(a^\varepsilon) \}_\varepsilon$ lies in a compact
subset of $W^{-1, p}(\Omega)$,

\vspace{0.1cm}

(iv) $\{\mbox{C}(b^\varepsilon)\}_\varepsilon$ lies in a compact
subset of $[W^{-1, q}(\Omega)]^{N \times N}$,

\vspace{0.1cm}

Then $a^\varepsilon b^\varepsilon \to ab$ in the sense of
distributions. }

\vspace{0.5cm}

{\bf Proof.} For any $\varepsilon  > 0 $, let $u^\varepsilon$ and
$v^\varepsilon$ be the weak solutions of the following problems:
$$
- \bigtriangleup u^\varepsilon = a^\varepsilon\quad \mbox{in}\;
\Omega , \;\quad u^\varepsilon = 0 \; \mbox{on} \; \partial \Omega,
\leqno (2.1)
$$
$$
- \bigtriangleup v^\varepsilon = b^\varepsilon\quad \mbox{in}\;
\Omega , \;\quad v^\varepsilon = 0 \; \mbox{on} \; \partial \Omega.
\leqno (2.2)
$$
From hypotheses \emph{(i)-(ii)} it follows that
$$
u^\varepsilon \rightharpoonup u \;\;\;\;\mbox{weakly in}\;\;\;
[W^{1, p} (\Omega)]^N,  \leqno(2.3)
$$
$$
v^\varepsilon \rightharpoonup v \;\;\;\;\mbox{weakly in}\;\;\;[W^{1,
q} (\Omega)]^N, \leqno(2.4)
$$
\noindent where $u$ and $v$ are uniquely defined by the problems:
$$
- \bigtriangleup u = a  \;\;\mbox{in}\;\; \Omega, \;\;\; u = 0
\;\;\mbox{on}\;\;
\partial \Omega, \leqno(2.5)
$$
$$
- \bigtriangleup v = b \;\;\mbox{in}\;\; \Omega, \;\;\; v = 0
\;\;\mbox{on}\;\;
\partial \Omega. \leqno(2.6)
$$

For $\varphi \in \mathcal{D}(\Omega),$ we denote by $\Omega_\varphi$
a smooth open set with the property:
$$
\mbox{supp} \;\varphi \subseteqq \Omega_\varphi \subseteqq \overline
\Omega_\varphi \subseteqq \Omega.\leqno (2.7)
$$

Then, the integral formula (1.1) yields:
$$
\begin{array}{l}
\int\limits_\Omega a^\varepsilon b^\varepsilon \varphi = \langle \;
\mbox{D}(a^\varepsilon), \varphi \; \mbox{D}(v^\varepsilon)\;
\rangle_p + \int\limits_{\Omega_\varphi} b^\varepsilon
\;\mbox{D}(u^\varepsilon)\bigtriangledown \varphi \;+ \\[7mm]
+ \int\limits_{\Omega_\varphi} \; \mbox{D}(u^\varepsilon)
\;(\bigtriangledown \varphi)\bigtriangledown \mbox{D}(v^\varepsilon)
- \int\limits_{\Omega_\varphi} \;
\mbox{D}(v^\varepsilon) \;(\bigtriangledown \varphi)\bigtriangledown \mbox{D}(u^\varepsilon) \; + \\[5mm]
+ \sum\limits_{1 \leq i, j \leq N} \; \int\limits_{\Omega_\varphi}
\;b^\varepsilon_i \;\mbox{C}_{ij}(u^\varepsilon)\;\dfrac{\partial
\varphi}{\partial x_j} + \dfrac{1}{2} \sum\limits_{1 \leq i, j \leq
N} \; \langle \;\mbox{C}_{ij}(b^\varepsilon), \varphi\;
\mbox{C}_{ij}(u^\varepsilon)\; \rangle_q .
\end{array} \leqno(2.8)
$$

\vspace{0.1cm}

From hypotheses \emph{(iii)-(iv)} we find that there exists a
subsequence (still denoted with $\varepsilon$) for which the
following properties hold:
$$
 \; \mbox{D}(a^\varepsilon) \to  \; \mbox{D}(a)\;\;\;\;
\mbox{strongly in} \;\;W^{-1, q}(\Omega),\leqno(2.9)
$$
$$
\; \mbox{C}(b^\varepsilon)  \to \; \mbox{C}(b)\;\;\;\;
\mbox{strongly in}\; \;[W^{-1, p}(\Omega)]^{N\times N}, \leqno(2.10)
$$

The theorem about the interior regularity of the solutions of the
problems (1)-(2)(see [7] \S 6.3) implies:
$$
\{u_\varepsilon\}_\varepsilon\; \mbox{ is bounded in}\;  [W^{2,
p}(\Omega_\varphi)]^N, \leqno (2.11)
$$
$$
\{v_\varepsilon\}_\varepsilon\; \mbox{ is bounded in}\;  [W^{2,
q}(\Omega_\varphi)]^N .\leqno (2.12)
$$

In the light of (2.9)-(2.10) and using the appropriate
Rellich-Kondrachov compactness results (see [7] \S 5.7) that follow
from (2.11)-(2.12), we see that every term of the right-hand side of
(2.8) represents a product of a strongly convergent (sub)sequence by
a weakly convergent one. Hence, the right-hand side of (2.8) is
converging to the corresponding formula which represents
$\int\limits_\Omega ab \varphi\;$. By the uniqueness of the limit,
the convergence holds on the entire sequence.

\vspace{1cm}

{\bf References.}

\vspace{0.1cm}

[1] Murat F., {\it Compacit\' e par compensation}, Ann. Scuola
Norm. Sup. Pisa Cl. Sci., {\bf 5(4)}, 1978.

[2] Murat F., {\it Compacit\' e par compensation II, Recent
Methods in Nonlinear Analysis} (De Giorgi E., Magenes E., Mosco U.
eds.), Pitagora Editriee, Bologna, 1979.

[3] Tartar L., {\it Compensated compactness and applications to
partial differential equations}, Heriot-Watt Sympos., vol. IV,
Pitman, New York, 1979.

[4] Dacorogna B., {\it Weak Continuity and Weak Lower
Semicontinuity of Nonlinear Functionals}, Springer-Verlag, Berlin,
1982.

[5] Tartar L., {\it The compensated compactness method applied to
systems of conservation lows, Systems of Nontinear PDE} (Ball J.
ed.), D. Reidel Pub. Co., Dordrecht, 1983.

[6] Evans L.C., {\it Weak Convergence Methods for  Nonlinear
Partial Differential Equations} (CBMS Regional Conference Series
in Mathematics, vol.74), American Mathematical Society,
Providence, 1991.

[7] Evans L.C., {\it Partial Differential Equations} (Graduate
Studies in Mathematics, vol.19), American Mathematical Society,
Providence, 1998.

\vspace{1cm}

\hspace{2.5cm}{\textbf{\dag\;\;}}\;{\sc I.M.A.R., P.O. Box 1-764,
Bucharest, Romania}

\end{document}